\documentclass[11pt]{article}
\usepackage{geometry}                
\usepackage{amsmath,amssymb}
\usepackage{amsthm}
\usepackage{hyperref}
\usepackage{graphicx,epsfig,color}
\usepackage{booktabs,bm,multirow}
\usepackage{cases}
\usepackage[caption=false]{subfig}

\newtheorem{theorem}{Theorem}[section]

\newtheorem{corollary}[theorem]{Corollary}

\newtheorem{remark}[theorem]{Remark}
\newtheorem{definition}[theorem]{Definition}
\newtheorem{lemma}[theorem]{Lemma}

\input{RREKTR.sty}

\title{Randomized regularized extended Kaczmarz algorithms for tensor recovery}
\author{Kui Du\thanks{School of Mathematical Sciences and Fujian Provincial Key Laboratory of Mathematical Modeling and High Performance Scientific Computing, Xiamen University, Xiamen 361005, China (kuidu@xmu.edu.cn).},
\quad  Xiao-Hui Sun\thanks{School of Mathematical Sciences, Xiamen University, Xiamen 361005, China ({sunxh@stu.xmu.edu.cn}).}} 
\date{}                                           

\begin{document}
\maketitle

\begin{abstract} Randomized regularized Kaczmarz algorithms have recently been proposed to solve tensor recovery models with {\it consistent} linear measurements. In this work, we propose a novel algorithm based on the randomized extended Kaczmarz algorithm (which converges linearly in expectation to the unique minimum norm least squares solution of a linear system) for tensor recovery models with {\it inconsistent} linear measurements. We prove the linear convergence in expectation of our algorithm. Numerical experiments on a tensor least squares problem and a sparse tensor recovery problem are given to illustrate the theoretical results.
\vspace{.5mm} 



\end{abstract}

\section{Introduction} Tensor recovery models have attracted much attention recently because of various applications, such as transportation, medical imaging, and remote sensing. With the assumption of {\it consistent} linear measurements, Chen and Qin \cite{chen2021regul} proposed a algorithmic framework based on the Kaczmarz-type algorithms \cite{kaczmarz1937angen,strohmer2009rando} for the following tensor recovery model: \beq\label{ctr} \wh\mcalx=\argmin_{\mcalx\in\mbbr^{N_2\times K\times N_3}}f(\mcalx),\quad \mbox{s.t.}\quad \mcala*\mcalx= \mcalb,\eeq where the objective function $f$ is strongly convex, the sensing tensor $\mcala\in\mbbr^{N_1\times N_2\times N_3}$, the acquired measurement tensor $\mcalb\in\mbbr^{N_1\times K\times N_3}$, and $\mcala*\mcalx$ is the tensor t-product \cite{kilmer2011facto}. In this paper, we assume that the linear measurement $\mcala*\mcalx= \mcalb$ is {\it inconsistent} and therefore consider the following constrained minimization problem: \beq\label{mp} \wh\mcalx=\argmin_{\mcalx\in\mbbr^{N_2\times K\times N_3}}f(\mcalx),\quad \mbox{s.t.}\quad \mcala^\top*\mcala*\mcalx=\mcala^\top*\mcalb,\eeq where $\mcala^\top$ denotes the transpose of $\mcala$ (see section 2.1). The solution $\wh\mcalx$ of (\ref{mp}) is a least squares solution of the inconsistent system $\mcala*\mcalx= \mcalb$ with some desirable characteristics promoted by regularization terms of the objective function $f$. 

In recent years, randomized iterative algorithms for linear systems of equations with massive data sets  have been greatly developed due to low memory footprints and good numerical performance, such as the randomized Kaczmarz (RK) algorithm \cite{strohmer2009rando}, the randomized coordinate descent (RCD) algorithm \cite{leventhal2010rando}, the randomized extended Kaczmarz (REK) algorithm \cite{zouzias2013rando}, and their extensions, e.g., \cite{needell2014paved,gower2015rando,ma2015conve,needell2015rando,bai2018greed,bai2019greed,bai2019parti,du2019tight,necoara2019faste,du2020rando,bai2021greed,wu2021two,zhang2021block,jiang2022rando,zhang2021greed}. The RK algorithm converges linearly in expectation to a solution of consistent linear systems \cite{strohmer2009rando,zouzias2013rando} and to within a radius ({\it convergence horizon}) of a (least squares) solution of inconsistent linear systems \cite{needell2010rando}. The REK algorithm converges linearly in expectation to a (least squares) solution of arbitrary linear systems \cite{zouzias2013rando,ma2015conve,du2019tight}. In this paper, we replace the RK algorithm integrated in the randomized regularized Kaczmarz (RRK) algorithm of \cite[Algorithm 3.2]{chen2021regul} with the REK algorithm and propose a randomized regularized extended Kaczmarz (RREK) algorithm for solving (\ref{mp}). The proposed RREK algorithm is called ``regularized'' since the objective function contains regularization terms for preserving some desirable characteristics of the underlying solution. We prove the linear convergence of the proposed algorithm. Special cases including tensor least squares problems and sparse tensor recovery problems are provided. Numerical experiments are given to illustrate our theoretical results.

The rest of this paper is organized as follows. In section 2, we provide clarification of notation, review basic concepts and results in tensor algebra and convex optimization, and also briefly introduce the RRK algorithm of \cite{chen2021regul}. In section 3, we describe the proposed RREK algorithm for solving (\ref{mp}) and establish its convergence theory. In section 4, we discuss two special cases including a tensor least squares problem and a sparse tensor recovery problem. In section 5, we report two numerical experiments to illustrate the theoretical results. Finally, we present brief concluding remarks in section 6.

\section{Preliminaries} 

\subsection{Basic notation}
Throughout the paper, we use 
boldface uppercase letters such as $\mbf A$ for matrices and calligraphic letters such as $\mcala$ for tensors. For an integer $m\geq 1$, let $[m]:=\{1,2,3,\ldots,m\}$. 
For any matrix $\mbf A\in\mbbr^{m\times n}$, we use $\mbf A^\top$, $\mbf A^\dag$, $\|\mbf A\|_2$, and $\sigma_{\rm min}(\mbf A)$ to denote the transpose, the Moore--Penrose pseudoinverse, the spectral norm, and the minimum nonzero singular values of $\mbf A$, respectively. For any random variables $\bm \xi$ and $\bm\zeta$, we use $\mbbe\bem\bm\xi\eem$ and $\mbbe\bem\bm\xi \mid \bm\zeta\eem$ to denote the expectation of $\bm\xi$ and the conditional expectation of $\bm\xi$ given $\bm\zeta$, respectively. 

\subsection{Tensor basics} In this subsection, we provide a brief review of key definitions and facts in tensor algebra. We follow the notation used in \cite{kilmer2011facto,kilmer2013third,miao2020gener}.

For a third-order tensor $\mcala\in\mbbr^{N_1\times N_2\times N_3}$, we denote its $(i,j,k)$ entry as $\mcala_{ijk}$ and use $\mcala_{i,:,:}$, $\mcala_{:,i,:}$ and $\mcala_{:,:,i}$ to denote respectively the $i$th horizontal, lateral and frontal slice. For notational convenience, the $i$th frontal slice $\mcala_{:,:,i}$ will be denoted as $\mbf A_i$. We define the block circulant matrix $\bcirc(\mcala)$ of $\mcala$ as, $$\bcirc(\mcala):=\bem\mbf A_1 &\mbf A_{N_3} & \cdots & \mbf A_2\\ \mbf A_2 & \mbf A_1 & \cdots &\mbf A_3\\ \vdots &\vdots & \ddots &\vdots\\ \mbf A_{N_3}& \mbf A_{N_3-1}& \cdots &\mbf A_1\eem\in\mbbr^{N_1N_3\times N_2N_3}.$$ We also define the operator $\unfold(\cdot)$ and its inversion $\fold(\cdot)$, $$\unfold(\mcala):=\bem\mbf A_1\\ \mbf A_2\\ \vdots\\ \mbf A_{N_3}\eem\in\mbbr^{N_1N_3\times N_2},\qquad \fold\l(\bem\mbf A_1\\ \mbf A_2\\ \vdots\\ \mbf A_{N_3}\eem\r):=\mcala.$$ 

\begin{definition}[t-product]
For $\mcala\in\mbbr^{N_1\times N_2\times N_3}$ and $\mcalb\in\mbbr^{N_2\times K\times N_3}$, the {\it t-product} $\mcala*\mcalb$ is defined to be a tensor of size $N_1\times K\times N_3$, $$\mcala*\mcalb:=\fold(\bcirc(\mcala)\unfold(\mcalb)).$$	
\end{definition}

\begin{definition}[transpose]
The transpose of $\mcala\in\mbbr^{N_1\times N_2\times N_3}$, denoted by $\mcala^\top$, is the $N_2\times N_1\times N_3$ tensor obtained by transposing each of the frontal slices and then reversing the order of transposed frontal slices $2$ through $N_3$.
\end{definition} 
For $\mcala\in\mbbr^{N_1\times N_2\times N_3}$, we have  \beq\label{ttop}\bcirc(\mcala^\top)=(\bcirc(\mcala))^\top.\eeq 
   
\begin{definition}[identity tensor]
 The identity tensor $\mcali\in\mbbr^{N\times N\times N_3}$ is the tensor whose first frontal slice is the $N\times N$ identity matrix, and whose other frontal slices are all zeros.	
\end{definition}
For $\mcala\in\mbbr^{N_1\times N_2\times N_3}$ and $\mcali\in\mbbr^{N_1\times N_1\times N_3}$, it holds that \beq\label{ii}\mcala_{i,:,:}=\mcali_{i,:,:}*\mcala.\eeq For $\mcala\in\mbbr^{N_1\times N_2\times N_3}$ and $\mcali\in\mbbr^{N_2\times N_2\times N_3}$, it holds that \beq\label{ij}\mcala_{:,j,:}=\mcala*\mcali_{:,j,:}.\eeq

\begin{definition}[inner product]
The inner product between $\mcala$ and $\mcalb$ in $\mbbr^{N_1\times N_2\times N_3}$ is defined as $$\la\mcala,\mcalb\ra:=\sum_{i,j,k}\mcala_{ijk}\mcalb_{ijk}.$$ 
\end{definition}

For $\mcala\in\mbbr^{N_1\times N_2\times N_3}$, $\mcalb\in\mbbr^{N_2\times K\times N_3}$, and $\mcalc\in\mbbr^{N_1\times K\times N_3}$, it holds that \beq\label{tran}\la\mcala*\mcalb,\mcalc\ra=\la\mcalb,\mcala^\top*\mcalc\ra.\eeq

\begin{definition}[$1$-norm, spectral norm, and Frobenius norm]
The 1-norm, spectral norm, and Frobenius norm of $\mcala\in\mbbr^{N_1\times N_2\times N_3}$ are defined as $$\|\mcala\|_1:=\sum_{i,j,k}|\mcala_{ijk}|,\qquad \|\mcala\|_2:=\|\bcirc(\mcala)\|_2,$$  and  $$\|\mcala\|_\rmf:=\sqrt{\la\mcala,\mcala\ra}=\sqrt{\sum_{i,j,k}(\mcala_{ijk})^2},$$ respectively.
\end{definition}

For $\mcala$ and $\mcalb$ in $\mbbr^{N_1\times N_2\times N_3}$, it holds that \beq\label{cs}|\la\mcala,\mcalb\ra|\leq \|\mcala\|_\rmf\|\mcalb\|_\rmf.\eeq For $\mcala\in\mbbr^{N_1\times N_2\times N_3}$ and $\mcalb\in\mbbr^{N_2\times K\times N_3}$, it holds that \beq\label{2f}\|\mcala*\mcalb\|_\rmf\leq\|\mcala\|_2\|\mcalb\|_\rmf.\eeq

\begin{definition}[$K$-range]
The $K$-range of $\mcala\in\mbbr^{N_1\times N_2\times N_3}$ is defined  as $$\ran_K(\mcala):=\{\mcala*\mcaly \mid \mcaly\in\mbbr^{N_2\times K\times N_3} \}.$$
\end{definition}
For $\mcala\in\mbbr^{N_1\times N_2\times N_3}$, it holds that $$\ran_K(\mcala^\top*\mcala)=\ran_K(\mcala^\top).$$ 
For $\mcala\in\mbbr^{N_1\times N_2\times N_3}$ and all $\mcalx\in \ran_K(\mcala)$, it holds that \beq\label{ran}\|\mcala^\top*\mcalx\|_\rmf^2\geq\sigma_{\min}^2(\bcirc(\mcala))\|\mcalx\|_\rmf^2.\eeq 


\begin{definition}[pseudoinverse] 
The pseudoinverse of $\mcala\in\mbbr^{N_1\times N_2\times N_3}$, denoted by $\mcala^\dag$, is the $N_2\times N_1\times N_3$ tensor  satisfying $$\bcirc(\mcala^\dag)=\bcirc(\mcala)^\dag.$$   	
\end{definition}
For $\mcala\in\mbbr^{N_1\times N_2\times N_3}$ and $\mcalb\in\mbbr^{N_1\times K\times N_3}$, it holds that $$\mcala^\top*\mcala*\mcala^\dag*\mcalb=\mcala^\top*\mcalb.$$ If $\mcalx\in\mbbr^{N_2\times K\times N_3}$ satisfies $\mcala^\top*\mcala*\mcalx=\mcala^\top*\mcalb$, then it holds that $$\mcala*\mcalx =\mcala*\mcala^\dag*\mcalb.$$  

\subsection{Convex optimization basics} To make the paper self-contained, we present basic definitions and properties about convex functions defined on tensor spaces in this subsection. We refer the reader to \cite{rockafellar1970conve,beck2017first} for more definitions and properties. 

\begin{definition}[subdifferential]
For a continuous function $f:\mbbr^{N_1\times N_2 \times N_3}\rightarrow\mbbr$, its subdifferential at $\mcalx\in\mbbr^{N_1\times N_2 \times N_3}$ is defined as $$\p f(\mcalx):=\{\mcalz:f(\mcaly)\geq f(\mcalx)+\la\mcalz,\mcaly-\mcalx\ra \quad  \forall\ \mcaly\in\mbbr^{N_1\times N_2 \times N_3}\}.$$	
\end{definition}

\begin{definition}[$\gamma$-strong convexity]
A function $f:\mbbr^{N_1\times N_2 \times N_3}\rightarrow\mbbr$ is called $\gamma$-strongly convex for a given $\gamma>0$ if the following inequality holds for all $\mcalx, \mcaly\in\mbbr^{N_1\times N_2\times N_3}$ and $\mcalz\in\p f(\mcalx)$: $$f(\mcaly)\geq f(\mcalx)+\la\mcalz,\mcaly-\mcalx\ra+\frac{\gamma}{2}\|\mcaly-\mcalx\|_\rmf^2.$$	
\end{definition}
The function $f(\mcalx)=\dsp\frac{1}{2}\|\mcalx\|_\rmf^2$ is differential and $1$-strongly convex. Moreover, it is easy to show that the function $h(\mcalx)+\dsp\frac{1}{2}\|\mcalx\|_\rmf^2$ is $1$-strongly convex if $h(\mcalx)$ is convex.


\begin{definition}[conjugate function]
The conjugate function of $f:\mbbr^{N_1\times N_2 \times N_3}\rightarrow\mbbr$ at $\mcaly\in\mbbr^{N_1\times N_2 \times N_3}$ is defined as $$f^*(\mcaly):=\sup_{\mcalx\in\mbbr^{N_1\times N_2 \times N_3}}\{\la\mcaly,\mcalx\ra-f(\mcalx)\}.$$	
\end{definition} 
If $f(\mcalx)$ is $\gamma$-strongly convex, then the conjugate function $f^*(\mcalx)$ is differentiable and for all $\mcalx,\mcaly\in\mbbr^{N_1\times N_2 \times N_3}$, the following inequality holds: 
 \beq\label{sc2}f^*(\mcaly)\leq f^*(\mcalx)+\la\nabla f^*(\mcalx),\mcaly-\mcalx\ra+\frac{1}{2\gamma}\|\mcaly-\mcalx\|_\rmf^2.\eeq For a strongly convex function $f(\mcalx)$, it can be shown that \cite{rockafellar1970conve,beck2017first} \beq\label{xgz}\mcalz\in\p f(\mcalx) \ \Leftrightarrow\ \mcalx=\nabla f^*(\mcalz).\eeq For a convex function $h(\mcalx)$, the conjugate function of $\lambda h(\mcalx)+\dsp\frac{1}{2}\|\mcalx\|_\rmf^2$ is differentiable. Its gradient involves the proximal mapping of $h(\mcalx)$, and it holds that $$\nabla\l(\lambda h(\cdot)+\frac{1}{2}\|\cdot\|_\rmf^2\r)^*(\mcalx)={\rm prox}_{\lambda h}(\mcalx):=\argmin_{\mcaly\in\mbbr^{N_1\times N_2\times N_3}}\l\{\lambda h(\mcaly)+\frac{1}{2}\|\mcalx-\mcaly\|_\rmf^2\r\}.$$

\begin{definition}[Bregman distance] For a convex function $f:\mbbr^{N_1\times N_2 \times N_3}\rightarrow\mbbr$, the Bregman distance between $\mcalx$ and $\mcaly$ with respect to $f$ and $\mcalz\in \p f(\mcalx)$ is defined as $$D_{f,\mcalz}(\mcalx,\mcaly):=f(\mcaly)-f(\mcalx)-\la\mcalz,\mcaly-\mcalx\ra.$$ 
\end{definition}
It follows from $\la\mcalz,\mcalx\ra=f(\mcalx)+f^*(\mcalz)$ if $\mcalz\in\p f(\mcalx)$ that \beq\label{dfz} D_{f,\mcalz}(\mcalx,\mcaly)=f(\mcaly)+f^*(\mcalz)-\la\mcalz,\mcaly\ra.\eeq If $f$ is $\gamma$-strongly convex, then it holds that \beq\label{gamma} D_{f,\mcalz}(\mcalx,\mcaly)\geq \frac{\gamma}{2}\|\mcalx-\mcaly\|_\rmf^2.\eeq


\begin{definition}[restricted strong convexity \cite{lai2013augme,schopfer2016linea}]
Let $f:\mbbr^{N_1\times N_2 \times N_3}\rightarrow\mbbr$ be convex differentiable with a nonempty minimizer set $X_f$. The function $f$ is called restricted strongly convex on $\mbbr^{N_1\times N_2\times N_3}$ with a constant $\mu>0$ if it satisfies for all $\mcalx\in \mbbr^{N_1\times N_2\times N_3}$ the inequality $$\la \nabla f(P_{X_f}(\mcalx))-\nabla f(\mcalx), P_{X_f}(\mcalx)-\mcalx\ra\geq \mu\|P_{X_f}(\mcalx)-\mcalx\|_\rmf^2,$$ where $P_{X_f}(\mcalx)$ denotes the orthogonal projection of $\mcalx$ onto $X_f$. 
\end{definition}

\begin{definition}[strong admissibility]
Let $\mcala\in\mbbr^{N_1\times N_2\times N_3}$ and $\mcalb\in\mbbr^{N_1\times K\times N_3}$ be given. Let $f:\mbbr^{N_2\times K \times N_3}\rightarrow\mbbr$ be strongly convex. The function $f$ is called strongly admissible if the function $g(\mcaly):=f^*(\mcala^\top*\mcala*\mcaly)-\la\mcaly,\mcala^\top*\mcalb\ra$ is restricted strongly convex on $\mbbr^{N_2\times K\times N_3}$.
\end{definition}
\begin{lemma}\label{lemmanu} Let $\wh\mcalx$ be the solution of (\ref{mp}). If $f$ is strongly admissible, then there exists a constant $\nu>0$ such that \beq\label{nu} D_{f,\mcalz}(\mcalx,\wh\mcalx)\leq \frac{1}{\nu}\|\mcala*(\mcalx-\wh\mcalx)\|_\rmf^2,\eeq for all $\mcalx\in \mbbr^{N_2\times K\times N_3}$ and $\mcalz\in\p f(\mcalx)\cap\ran_K(\mcala^\top)$.
\end{lemma}	
\proof
 The solution $\wh\mcalx$ of (\ref{mp}) satisfies the following optimality conditions: \beq\label{oc}\mcala^\top*\mcala*\wh\mcalx=\mcala^\top*\mcalb,\qquad \p f(\wh\mcalx)\cap\ran_K(\mcala^\top*\mcala)\neq\emptyset.\eeq The dual problem of (\ref{mp}) is the unconstrained problem $$\min_{\mcaly\in\mbbr^{N_2\times K\times N_3}}g(\mcaly),$$ where $$ g(\mcaly)=f^*(\mcala^\top*\mcala*\mcaly)-\la\mcaly,\mcala^\top*\mcalb\ra.$$ By the strong duality, we have $$f(\wh\mcalx)=-\min_{\mcaly \in\mbbr^{N_2\times K\times N_3}}g(\mcaly).$$ Since $\mcalz\in\ran_K(\mcala^\top)=\ran_K(\mcala^\top*\mcala)$, we can write $\mcalz=\mcala^\top*\mcala*\mcaly$ for some $\mcaly$. Then \begin{align*} D_{f,\mcalz}(\mcalx,\wh\mcalx) &  \stackrel{(\ref{dfz})}{=} f^*(\mcalz)-\la\mcalz,\wh\mcalx\ra+f(\wh\mcalx)\\ &=f^*(\mcala^\top*\mcala*\mcaly) -\la\mcala^\top*\mcala*\mcaly,\wh\mcalx\ra +f(\wh\mcalx)\\ &\stackrel{(\ref{tran})}= f^*(\mcala^\top*\mcala*\mcaly) -\la\mcaly,\mcala^\top*\mcala*\wh\mcalx\ra +f(\wh\mcalx)\\ &\stackrel{(\ref{oc})}{=} f^*(\mcala^\top*\mcala*\mcaly)-\la\mcaly,\mcala^\top*\mcalb\ra +f(\wh\mcalx)\\ &= g(\mcaly)-\min_{\mcaly \in\mbbr^{N_2\times K\times N_3}}g(\mcaly).\end{align*} Since $g(\mcaly)$ is restricted strongly convex on $\mbbr^{N_2\times K\times N_3}$, there exists a constant $\mu>0$ such that $$\la \nabla g(P_{X_g}(\mcaly))-\nabla g(\mcaly), P_{X_g}(\mcaly)-\mcaly\ra\geq \mu\|P_{X_g}(\mcaly)-\mcaly\|_\rmf^2.$$ By $\nabla g(P_{X_g}(\mcaly))=0$ and the Cauchy--Schwarz inequality (\ref{cs}), we get $$\|\nabla g(\mcaly)\|_\rmf\geq\mu\|P_{X_g}(\mcaly)-\mcaly\|_\rmf.$$  The convexity of $g(\mcaly)$ implies \beq\label{convex} g(\mcaly)-g(P_{X_g}(\mcaly))\leq \la\nabla g(\mcaly),\mcaly-P_{X_g}(\mcaly)\ra\leq\|\nabla g(\mcaly)\|_\rmf\|P_{X_g}(\mcaly)-\mcaly\|_\rmf\leq \frac{1}{\mu}\|\nabla g(\mcaly)\|_\rmf^2.\eeq The gradient of $g(\mcaly)$ is \begin{align} \nabla g(\mcaly)&=\mcala^\top*\mcala*\nabla f^*(\mcala^\top*\mcala*\mcaly)-\mcala^\top*\mcalb \nn \\ &= \mcala^\top*\mcala*\nabla f^*(\mcalz)-\mcala^\top*\mcalb \nn \\ & \stackrel{(\ref{xgz})}{=} \mcala^\top*\mcala*\mcalx-\mcala^\top*\mcalb \nn \\ &= \mcala^\top*\mcala*\mcalx-\mcala^\top*\mcala*\wh\mcalx. \label{grad}\end{align} Therefore, \begin{align*} D_{f,\mcalz}(\mcalx,\wh\mcalx) & =g(\mcaly)-g(P_{X_g}(\mcaly))\stackrel{(\ref{convex})(\ref{grad})}\leq \frac{1}{\mu}\|\mcala^\top*\mcala*\mcalx-\mcala^\top*\mcala*\wh\mcalx\|_\rmf^2\\ &\stackrel{(\ref{2f})}\leq \frac{\|\mcala\|_2^2}{\mu}\|\mcala*(\mcalx-\wh\mcalx)\|_\rmf^2:=\frac{1}{\nu}\|\mcala*(\mcalx-\wh\mcalx)\|_\rmf^2.\end{align*} This completes the proof. \endproof

 The constant $\nu$ depends on the tensor $\mcala$ and the function $f$. In general, it is hard to quantify $\nu$. We refer the reader to \cite{schopfer2016linea,lai2013augme, chen2021regul} for examples of strongly admissible functions.

\subsection{The RRK algorithm}
By combining the RK algorithm and the gradient of the conjugate function at the previous iterate, Chen and Qin \cite{chen2021regul} proposed the RRK algorithm (see Algorithm 1) for solving the minimization problem (\ref{ctr}). They proved a linear convergence rate if $\mcala*\mcalx=\mcalb$ is consistent; see Theorem 3.9 of \cite{chen2021regul}. Moreover, they also considered the noisy scenario (the perturbed constraint $\mcala*\mcalx=\wt\mcalb$ where $\wt\mcalb=\mcalb+\mcale$) and proved that the RRK algorithm linearly converges to with a radius of the solution of (\ref{ctr}); see Theorem 3.10 of \cite{chen2021regul}. 

\begin{center}
\begin{tabular*}{160mm}{l}
\toprule {\bf Algorithm 1:} The RRK algorithm for solving (\ref{ctr})\\ 
\hline \noalign{\smallskip}
\quad {\bf Input}: $\mcala\in\mbbr^{N_1\times N_2\times N_3}$, $\mcalb\in\mbbr^{N_1\times K\times N_3}$, stepsize $\alpha_{\rm r}$, maximum number of iterations {\tt M}, \\ \quad and tolerance $\tau$.\\\noalign{\smallskip}
\quad {\bf Initialize}: $\mcalz^{(0)}\in\ran_K(\mcala^\top)$ and $\mcalx^{(0)}=\nabla f^*(\mcalz^{(0)})$.\\ \noalign{\smallskip}
\quad {\bf for} $k=1,2,\ldots,$ {\tt M} {\bf do}\\ \noalign{\smallskip}
\quad \qquad  Pick $i_k\in[N_1]$ with probability ${\|\mcala_{i_k,:,:}\|^2_\rmf}/{\|\mcala\|_\rmf^2}$\\  \noalign{\smallskip}
\quad \qquad Set $\mcalz^{(k)} = \mcalz^{(k-1)}-\alpha_{\rm r}(\mcala_{i_k,:,:})^\top*\dsp \frac{\mcala_{i_k,:,:}*\mcalx^{(k-1)}-\mcalb_{i_k,:,:}}{\|\mcala_{i_k,:,:}\|^2_\rmf}$ \\  \noalign{\smallskip}
\quad \qquad Set $\mcalx^{(k)}=\nabla f^*(\mcalz^{(k)})$ \\  \noalign{\smallskip}
\quad \qquad Stop if $\|\mcalx^{(k)}-\mcalx^{(k-1)}\|_\rmf/\|\mcalx^{(k-1)}\|_\rmf<\tau$ \\  \noalign{\smallskip}
\quad {\bf end}\\
\bottomrule
\end{tabular*}
\end{center}

For the special case that $$f(\mcalx)=\frac{1}{2}\|\mcalx\|_\rmf^2,$$ we have $$f^*(\mcalx)=\frac{1}{2}\|\mcalx\|_\rmf^2,\qquad \nabla f^*(\mcalx)=\mcalx.$$ Algorithm 1 becomes a tensor randomized Kaczmarz (TRK) algorithm for solving the tensor system $\mcala*\mcalx=\mcalb$ with the iteration \beq\label{trk}\mcalx^{(k)} = \mcalx^{(k-1)}-\alpha_{\rm r}(\mcala_{i_k,:,:})^\top*\dsp \frac{\mcala_{i_k,:,:}*\mcalx^{(k-1)}-\mcalb_{i_k,:,:}}{\|\mcala_{i_k,:,:}\|^2_\rmf},\eeq which is different from the tensor randomized Kaczmarz algorithm proposed by Ma and Molitor \cite{ma2021rando}.

\section{The RREK algorithm}

The REK algorithm \cite{zouzias2013rando} converges linearly in expectation to the minimum 2-norm least squares solution of a linear system. Motivated by this property, we consider replacing the RK algorithm integrated in the RRK algorithm with the REK algorithm and propose the following RREK algorithm (see Algorithm 2) for solving the constrained minimization problem (\ref{mp}). The RREK algorithm only uses one horizontal slice and one lateral slice of $\mcala$ at each step and avoids forming $\mcala^\top*\mcala$ explicitly. We also note that $\mcalz^{(k)}$ in the RREK algorithm is the same as the $k$th iterate generated by the TRK iteration (\ref{trk}) applied to the consistent tensor system $\mcala^\top * \mcalz=\mbf 0$. 

\begin{center}
\begin{tabular*}{160mm}{l}
\toprule {\bf Algorithm 2:} The RREK algorithm for solving (\ref{mp})\\ 
\hline \noalign{\smallskip}
\quad {\bf Input}: $\mcala\in\mbbr^{N_1\times N_2\times N_3}$, $\mcalb\in\mbbr^{N_1\times K\times N_3}$, stepsizes $\alpha_{\rm r}$ and $\alpha_{\rm c}$, maximum number of \\ \quad iterations {\tt M}, and tolerance $\tau$.\\\noalign{\smallskip}
\quad {\bf Initialize}: $\mcalz^{(0)}=\mcalb$, $\mcaly^{(0)}\in\ran_K(\mcala^\top)$, and $\mcalx^{(0)}=\nabla f^*(\mcaly^{(0)})$.\\ \noalign{\smallskip}
\quad {\bf for} $k=1,2,\ldots,$ {\tt M} {\bf do}\\ \noalign{\smallskip}
\quad \qquad  Pick $j_k\in[N_2]$ with probability ${\|\mcala_{:,j_k,:}\|^2_\rmf}/{\|\mcala\|_\rmf^2}$\\  \noalign{\smallskip}
\quad \qquad  Set $\mcalz^{(k)}=\mcalz^{(k-1)}-\alpha_{\rm c}\mcala_{:,j_k,:}*\dsp\frac{(\mcala_{:,j_k,:})^\top*\mcalz^{(k-1)}}{\|\mcala_{:,j_k,:}\|_\rmf^2}$\\  \noalign{\smallskip}
\quad \qquad  Pick $i_k\in[N_1]$ with probability ${\|\mcala_{i_k,:,:}\|^2_\rmf}/{\|\mcala\|_\rmf^2}$\\  \noalign{\smallskip}
\quad \qquad Set $\mcaly^{(k)} = \mcaly^{(k-1)}-\alpha_{\rm r}(\mcala_{i_k,:,:})^\top*\dsp \frac{\mcala_{i_k,:,:}*\mcalx^{(k-1)}-\mcalb_{i_k,:,:}+\mcalz^{(k)}_{i_k,:,:}}{\|\mcala_{i_k,:,:}\|^2_\rmf}$ \\  \noalign{\smallskip}
\quad \qquad Set $\mcalx^{(k)}=\nabla f^*(\mcaly^{(k)})$ \\  \noalign{\smallskip}
\quad \qquad Stop if $\|\mcalx^{(k)}-\mcalx^{(k-1)}\|_\rmf/\|\mcalx^{(k-1)}\|_\rmf<\tau$ \\  \noalign{\smallskip}
\quad {\bf end}\\
\bottomrule
\end{tabular*}
\end{center}

Next we analyze the convergence of the RREK algorithm. Our analysis is similar to that of \cite{du2020rando}, but slightly more complicated. The convergence estimates depend on the positive numbers $\lambda_{\rm r}$ and $\lambda_{\rm c}$ defined as $$\lambda_{\rm r}:=\max_{i\in[N_1]}\frac{\|\mcala_{i,:,:}\|_2^2}{\|\mcala_{i,:,:}\|_\rmf^2},\qquad \lambda_{\rm c}:=\max_{j\in[N_2]}\frac{\|\mcala_{:,j,:}\|_2^2}{\|\mcala_{:,j,:}\|_\rmf^2}.$$ We give the convergence result of $\mcalz^{(k)}$ in the RREK algorithm in the following theorem.

\begin{theorem}\label{atz} If $0<\alpha_{\rm c}<2/\lambda_{\rm c}$, then the sequence $\{\mcalz^{(k)}\}$ generated by Algorithm 2 satisfies $$\mbbe\bem\|\mcalz^{(k)}-(\mcalb-\mcala*\mcala^\dag*\mcalb)\|_\rmf^2\eem\leq \rho_{\rm c}^k \|\mcala*\mcala^\dag*\mcalb\|_\rmf^2,$$ where $$\rho_{\rm c}=1-\frac{(2\alpha_{\rm c}-\alpha_{\rm c}^2\lambda_{\rm c})\sigma_{\min}^2(\bcirc(\mcala))}{\|\mcala\|_\rmf^2}.$$
\end{theorem}

\proof
Introduce the auxiliary tensor sequence $$\mcale^{(k)}=\mcalz^{(k)}-(\mcalb-\mcala*\mcala^\dag*\mcalb).$$ By $(\mcala_{:,j_k,:})^\top*(\mcalb-\mcala*\mcala^\dag*\mcalb)=\mbf 0$ (because $\mcala^\top*(\mcalb-\mcala*\mcala^\dag*\mcalb)=\mbf 0$), we have $$\mcale^{(k)} = \mcale^{(k-1)}-\alpha_{\rm c}\mcala_{:,j_k,:}*\dsp\frac{(\mcala_{:,j_k,:})^\top*\mcale^{(k-1)}}{\|\mcala_{:,j_k,:}\|_\rmf^2}.$$ Then, \begin{align*}
 \|\mcale^{(k)}\|_\rmf^2 & =\l\|\mcale^{(k-1)}-\alpha_{\rm c}\mcala_{:,j_k,:}*\dsp\frac{(\mcala_{:,j_k,:})^\top*\mcale^{(k-1)}}{\|\mcala_{:,j_k,:}\|_\rmf^2}\r\|_\rmf^2\\ &=\|\mcale^{(k-1)}\|_\rmf^2-\frac{2\alpha_{\rm c}}{\|\mcala_{:,j_k,:}\|_\rmf^2}\la\mcale^{(k-1)},\mcala_{:,j_k,:}*(\mcala_{:,j_k,:})^\top*\mcale^{(k-1)}\ra	\\ &\quad\ +\frac{\alpha_{\rm c}^2}{\|\mcala_{:,j_k,:}\|_\rmf^4}\|\mcala_{:,j_k,:}*(\mcala_{:,j_k,:})^\top*\mcale^{(k-1)}\|_\rmf^2\\ &\stackrel{(\ref{tran})(\ref{2f})}\leq \|\mcale^{(k-1)}\|_\rmf^2-\frac{2\alpha_{\rm c}}{\|\mcala_{:,j_k,:}\|_\rmf^2}\|(\mcala_{:,j_k,:})^\top*\mcale^{(k-1)}\|_\rmf^2+\frac{\alpha_{\rm c}^2\|\mcala_{:,j_k,:}\|_2^2}{\|\mcala_{:,j_k,:}\|_\rmf^4}\|(\mcala_{:,j_k,:})^\top*\mcale^{(k-1)}\|_\rmf^2\\ 
 &\leq \|\mcale^{(k-1)}\|_\rmf^2-\frac{2\alpha_{\rm c}}{\|\mcala_{:,j_k,:}\|_\rmf^2}\|(\mcala_{:,j_k,:})^\top*\mcale^{(k-1)}\|_\rmf^2 +\frac{\alpha_{\rm c}^2\lambda_{\rm c}}{\|\mcala_{:,j_k,:}\|_\rmf^2}\|(\mcala_{:,j_k,:})^\top*\mcale^{(k-1)}\|_\rmf^2.
 \end{align*}
 Taking conditional expectation conditioned on $\mcale^{(k-1)}$ gives
 \begin{align*}
 \mbbe\bem\|\mcale^{(k)}\|_\rmf^2\mid \mcale^{(k-1)}\eem &\stackrel{(\ref{ij})}\leq \|\mcale^{(k-1)}\|_\rmf^2-\frac{2\alpha_{\rm c}-\alpha_{\rm c}^2\lambda_{\rm c}}{\|\mcala\|_\rmf^2}\|\mcala^\top*\mcale^{(k-1)}\|_\rmf^2\\ & \leq \l(1-\frac{(2\alpha_{\rm c}-\alpha_{\rm c}^2\lambda_{\rm c})\sigma_{\min}^2(\bcirc(\mcala))}{\|\mcala\|_\rmf^2}\r)\|\mcale^{(k-1)}\|_\rmf^2.	
 \end{align*} In the last inequality, we use the facts that $2\alpha_{\rm c}-\alpha_{\rm c}^2\lambda_{\rm c}>0$, $\mcale^{(k)}\in\ran_K(\mcala)$ (by induction), and (\ref{ran}). Next, by the law of total expectation, we have
\begin{align*}
 \mbbe\bem\|\mcale^{(k)}\|_\rmf^2\eem \leq \l(1-\frac{(2\alpha_{\rm c}-\alpha_{\rm c}^2\lambda_{\rm c})\sigma_{\min}^2(\bcirc(\mcala))}{\|\mcala\|_\rmf^2}\r)\mbbe\bem\|\mcale^{(k-1)}\|_\rmf^2\eem.	
\end{align*}
Unrolling the recurrence yields the result.
\endproof

We give the main convergence result of the RREK algorithm in the following theorem.
\begin{theorem}\label{main} Let $f$ be $\gamma$-strongly convex and strongly admissible, and $\nu>0$ is the constant  from (\ref{nu}). Let $\wh\mcalx$ be the solution of (\ref{mp}). Assume that $0<\alpha_{\rm c}<2/\lambda_{\rm c}$ and $0<\alpha_{\rm r}<2\gamma/\lambda_{\rm r}$. For any $\delta>0$, the sequences $\{\mcalx^{(k)}\}$ and $\{\mcaly^{(k)}\}$ generated by Algorithm 2 satisfy \begin{align*}\mbbe\bem D_{f,\mcaly^{(k)}}(\mcalx^{(k)},\wh\mcalx)\eem &\leq \frac{\delta+\gamma}{2\delta\gamma}\frac{\alpha_{\rm r}^2\lambda_{\rm r}}{\|\mcala\|_\rmf^2}\|\mcala*\mcala^\dag*\mcalb\|_\rmf^2\sum_{i=0}^{k-1}\rho_{\rm c}^{k-i}\l(1+\frac{\delta}{\gamma}\r)^i\rho_{\rm r}^i\\ &\quad +\l(1+\frac{\delta}{\gamma}\r)^k\rho_{\rm r}^kD_{f,\mcaly^{(0)}}(\mcalx^{(0)},\wh\mcalx),\end{align*} where $$\rho_{\rm c}=1-\frac{(2\alpha_{\rm c}-\alpha_{\rm c}^2\lambda_{\rm c})\sigma_{\min}^2(\bcirc(\mcala))}{\|\mcala\|_\rmf^2},\quad \rho_{\rm r}=1-\frac{(2\gamma\alpha_{\rm r}-\alpha_{\rm r}^2\lambda_{\rm r})\nu}{2\gamma\|\mcala\|_\rmf^2}.$$\end{theorem}

\proof
Let $$\wh\mcaly^{(k)} := \mcaly^{(k-1)}-\alpha_{\rm r}(\mcala_{i_k,:,:})^\top* \frac{\mcala_{i_k,:,:}*\mcalx^{(k-1)}-\mcala_{i_k,:,:}*\mcala^\dag*\mcalb}{\|\mcala_{i_k,:,:}\|^2_\rmf},$$ which is actually the one-step RRK update for the consistent constraint $$\mcala*\mcaly = \mcala*\mcala^\dag*\mcalb$$ from $\mcalx^{(k-1)}$ and $\mcaly^{(k-1)}.$  We have $$\mcaly^{(k)}-\wh\mcaly^{(k)}=\alpha_{\rm r}(\mcala_{i_k,:,:})^\top*\frac{\mcalb_{i_k,:,:}-\mcalz^{(k)}_{i_k,:,:}-\mcala_{i_k,:,:}*\mcala^\dag*\mcalb}{\|\mcala_{i_k,:,:}\|^2_\rmf}.$$
Then, \begin{align}
 	\|\mcaly^{(k)}-\wh\mcaly^{(k)}\|_\rmf^2 & = \frac{\alpha_{\rm r}^2}{\|\mcala_{i_k,:,:}\|_\rmf^4}\|(\mcala_{i_k,:,:})^\top*(\mcalb_{i_k,:,:}-\mcalz^{(k)}_{i_k,:,:}-\mcala_{i_k,:,:}*\mcala^\dag*\mcalb)\|_\rmf^2 \nn \\ &\stackrel{(\ref{ttop})(\ref{2f})}\leq \frac{\alpha_{\rm r}^2\|\mcala_{i_k,:,:}\|_2^2}{\|\mcala_{i_k,:,:}\|_\rmf^4}\|\mcalb_{i_k,:,:}-\mcalz^{(k)}_{i_k,:,:}-\mcala_{i_k,:,:}*\mcala^\dag*\mcalb\|_\rmf^2\nn \\ &\leq \frac{\alpha_{\rm r}^2\lambda_{\rm r}}{\|\mcala_{i_k,:,:}\|_\rmf^2}\|\mcalb_{i_k,:,:}-\mcalz^{(k)}_{i_k,:,:}-\mcala_{i_k,:,:}*\mcala^\dag*\mcalb\|_\rmf^2.
\label{yy} \end{align} Let $\mbbe_{k-1}\bem\cdot\eem$ denote the conditional expectation conditioned on $\mbf \mcalz^{(k-1)}$, $\mcaly^{(k-1)}$, and $\mcalx^{(k-1)}$. Let $\mbbe_{k-1}^{\rm r}\bem\cdot\eem$ denote the conditional expectation conditioned on $\mcalz^{(k)}$, $\mcaly^{(k-1)}$, and $\mcalx^{(k-1)}$. Then, by the law of total expectation, we have $$\mbbe_{k-1}\bem\cdot\eem=\mbbe_{k-1}\bem\mbbe_{k-1}^{\rm r}\bem\cdot\eem\eem.$$
 Taking conditional expectation for (\ref{yy}) conditioned on $\mcalz^{(k-1)}$, $\mcaly^{(k-1)}$, and $\mcalx^{(k-1)}$, we obtain
\begin{align*}\mbbe_{k-1}\bem\|\mcaly^{(k)}-\wh\mcaly^{(k)}\|_\rmf^2\eem &=\mbbe_{k-1}\bem\mbbe_{k-1}^{\rm r}\bem\|\mcaly^{(k)}-\wh\mcaly^{(k)}\|_\rmf^2\eem\eem\\ & \stackrel{(\ref{ii})(\ref{yy})}\leq \frac{\alpha_{\rm r}^2\lambda_{\rm r}}{\|\mcala\|_\rmf^2}\mbbe_{k-1}\bem\|\mcalz^{(k)}-(\mcalb-\mcala*\mcala^\dag*\mcalb)\|_\rmf^2 \eem. 
\end{align*}
Then, by the law of total expectation and Theorem \ref{atz}, we have
\begin{align}\mbbe\bem\|\mcaly^{(k)}-\wh\mcaly^{(k)}\|_\rmf^2\eem &\leq \frac{\alpha_{\rm r}^2\lambda_{\rm r}}{\|\mcala\|_\rmf^2}\mbbe\bem\|\mcalz^{(k)}-(\mcalb-\mcala*\mcala^\dag*\mcalb)\|_\rmf^2 \eem \nn\\ & \leq \frac{\alpha_{\rm r}^2\lambda_{\rm r}\rho_{\rm c}^k}{\|\mcala\|_\rmf^2}\|\mcala*\mcala^\dag*\mcalb\|_\rmf^2.\label{ey}\end{align}
Let 
$$\mcalw:= \frac{\mcala_{i_k,:,:}*\mcalx^{(k-1)}-\mcala_{i_k,:,:}*\mcala^\dag*\mcalb}{\|\mcala_{i_k,:,:}\|^2_\rmf}.$$ By $\mcala*\wh\mcalx=\mcala*\mcala^\dag*\mcalb$, we have \beq\label{w}\mcalw=\frac{\mcala_{i_k,:,:}*(\mcalx^{(k-1)}-\wh\mcalx)}{\|\mcala_{i_k,:,:}\|^2_\rmf}.\eeq
Let \beq\label{xhat}\wh\mcalx^{(k)}:=\nabla f^*(\wh\mcaly^{(k)}).\eeq By (\ref{xgz}), we have \beq\label{yhat}\wh\mcaly^{(k)}\in\p f(\wh\mcalx^{(k)}).\eeq 
Therefore, the Bregman distance between $\wh\mcalx^{(k)}$ and $\wh\mcalx$ with respect to $f$ and $\wh\mcaly^{(k)}$ satisfies \begin{align*}
D_{f,\wh\mcaly^{(k)}}(\wh\mcalx^{(k)},\wh\mcalx) &\stackrel{(\ref{dfz})}= f^*(\wh\mcaly^{(k)})-\la\wh\mcaly^{(k)},\wh\mcalx\ra+f(\wh\mcalx) \\ & = f^*( \mcaly^{(k-1)}-\alpha_{\rm r}(\mcala_{i_k,:,:})^\top*\mcalw)	-\la \mcaly^{(k-1)}-\alpha_{\rm r}(\mcala_{i_k,:,:})^\top*\mcalw,\wh\mcalx\ra+f(\wh\mcalx) \\ & \stackrel{(\ref{sc2})}\leq f^*(\mcaly^{(k-1)}) -\la\nabla f^*(\mcaly^{(k-1)}),\alpha_{\rm r}(\mcala_{i_k,:,:})^\top*\mcalw \ra+\frac{\alpha_{\rm r}^2}{2\gamma}\|(\mcala_{i_k,:,:})^\top*\mcalw\|_\rmf^2\\ &\qquad	-\la \mcaly^{(k-1)}-\alpha_{\rm r}(\mcala_{i_k,:,:})^\top*\mcalw,\wh\mcalx\ra+f(\wh\mcalx) \\ & \stackrel{(\ref{dfz})}= D_{f,\mcaly^{(k-1)}}(\mcalx^{(k-1)},\wh\mcalx)-\la\mcalx^{(k-1)},\alpha_{\rm r}(\mcala_{i_k,:,:})^\top*\mcalw\ra\\ &\qquad	+\frac{\alpha_{\rm r}^2}{2\gamma}\|(\mcala_{i_k,:,:})^\top*\mcalw\|_\rmf^2+\la\alpha_{\rm r}(\mcala_{i_k,:,:})^\top*\mcalw,\wh\mcalx\ra\\ & =  D_{f,\mcaly^{(k-1)}}(\mcalx^{(k-1)},\wh\mcalx)-\la\mcalx^{(k-1)}-\wh\mcalx,\alpha_{\rm r}(\mcala_{i_k,:,:})^\top*\mcalw\ra\\ & \qquad +\frac{\alpha_{\rm r}^2}{2\gamma}\|(\mcala_{i_k,:,:})^\top*\mcalw\|_\rmf^2\\ & \stackrel{(\ref{w})(\ref{tran})}  = D_{f,\mcaly^{(k-1)}}(\mcalx^{(k-1)},\wh\mcalx)-\frac{\alpha_{\rm r}}{\|\mcala_{i_k,:,:}\|_\rmf^2}\|\mcala_{i_k,:,:}*(\mcalx^{(k-1)}-\wh\mcalx)\|_\rmf^2 \\ &\qquad +\frac{\alpha_{\rm r}^2}{2\gamma}\|(\mcala_{i_k,:,:})^\top*\mcalw\|_\rmf^2 \\ & \stackrel{(\ref{ttop})(\ref{2f})} \leq D_{f,\mcaly^{(k-1)}}(\mcalx^{(k-1)},\wh\mcalx)-\frac{\alpha_{\rm r}}{\|\mcala_{i_k,:,:}\|_\rmf^2}\|\mcala_{i_k,:,:}*(\mcalx^{(k-1)}-\wh\mcalx)\|_\rmf^2 \\ &\qquad +\frac{\alpha_{\rm r}^2\lambda_{\rm r}}{2\gamma\|\mcala_{i_k,:,:}\|_\rmf^2}\|\mcala_{i_k,:,:}*(\mcalx^{(k-1)}-\wh\mcalx)\|_\rmf^2.
\end{align*} 
Taking conditional expectation conditioned on $\mcalz^{(k-1)}$, $\mcaly^{(k-1)}$, and $\mcalx^{(k-1)}$, we have 
\begin{align*}
\mbbe_{k-1}\bem D_{f,\wh\mcaly^{(k)}}(\wh\mcalx^{(k)},\wh\mcalx) \eem &\stackrel{(\ref{ii})}\leq D_{f,\mcaly^{(k-1)}}(\mcalx^{(k-1)},\wh\mcalx) -\frac{2\gamma\alpha_{\rm r}-\alpha_{\rm r}^2\lambda_{\rm r}}{2\gamma\|\mcala\|_\rmf^2}\|\mcala*(\mcalx^{(k-1)}-\wh\mcalx)\|_\rmf^2\\ &\stackrel{(\ref{xgz})(\ref{nu})}\leq \l(1-\frac{(2\gamma\alpha_{\rm r}-\alpha_{\rm r}^2\lambda_{\rm r})\nu}{2\gamma\|\mcala\|_\rmf^2}\r)D_{f,\mcaly^{(k-1)}}(\mcalx^{(k-1)},\wh\mcalx).
\end{align*}
Thus,  by the law of total expectation, we have
\beq\label{ed} \mbbe \bem D_{f,\wh\mcaly^{(k)}}(\wh\mcalx^{(k)},\wh\mcalx)\eem\leq \rho_{\rm r}\mbbe\bem D_{f,\mcaly^{(k-1)}}(\mcalx^{(k-1)},\wh\mcalx)\eem.\eeq
Now, we consider the Bregman distance $D_{f,\mcaly^{(k)}}(\mcalx^{(k)},\wh\mcalx)$, which satisfies \begin{align*} D_{f,\mcaly^{(k)}}(\mcalx^{(k)},\wh\mcalx) &\stackrel{(\ref{dfz})} =D_{f,\wh\mcaly^{(k)}}(\wh\mcalx^{(k)},\wh\mcalx)+f^*(\mcaly^{(k)})-f^*(\wh\mcaly^{(k)})-\la\mcaly^{(k)},\wh\mcalx\ra+\la\wh\mcaly^{(k)},\wh\mcalx\ra\\
& \stackrel{(\ref{sc2})}\leq D_{f,\wh\mcaly^{(k)}}(\wh\mcalx^{(k)},\wh\mcalx)+\la\nabla f^*(\wh\mcaly^{(k)}),\mcaly^{(k)}-\wh\mcaly^{(k)}\ra+\frac{1}{2\gamma}\|\mcaly^{(k)}-\wh\mcaly^{(k)}\|_\rmf^2\\ &\quad -\la\mcaly^{(k)}-\wh\mcaly^{(k)},\wh\mcalx\ra
\\ & \stackrel{(\ref{xhat})}= D_{f,\wh\mcaly^{(k)}}(\wh\mcalx^{(k)},\wh\mcalx)+\la\wh\mcalx^{(k)}-\wh\mcalx,\mcaly^{(k)}-\wh\mcaly^{(k)}\ra+\frac{1}{2\gamma}\|\mcaly^{(k)}-\wh\mcaly^{(k)}\|_\rmf^2
\\ & \leq  D_{f,\wh\mcaly^{(k)}}(\wh\mcalx^{(k)},\wh\mcalx)+\frac{\delta}{2}\|\wh\mcalx^{(k)}-\wh\mcalx\|_\rmf^2+\frac{1}{2\delta}\|\mcaly^{(k)}-\wh\mcaly^{(k)}\|_\rmf^2\\ &\qquad +\frac{1}{2\gamma}\|\mcaly^{(k)}-\wh\mcaly^{(k)}\|_\rmf^2 \quad \mbox{(by  (\ref{cs}) and Young's inequality) }
\\ &\stackrel{(\ref{yhat})(\ref{gamma})}\leq \l(1+\frac{\delta}{\gamma}\r)D_{f,\wh\mcaly^{(k)}}(\wh\mcalx^{(k)},\wh\mcalx)+\frac{\delta+\gamma}{2\delta\gamma}\|\mcaly^{(k)}-\wh\mcaly^{(k)}\|_\rmf^2.
\end{align*}
Taking expectation, we have
\begin{align*}
	\mbbe\bem D_{f,\mcaly^{(k)}}(\mcalx^{(k)},\wh\mcalx)\eem & \leq  \l(1+\frac{\delta}{\gamma}\r)\mbbe\bem D_{f,\wh\mcaly^{(k)}}(\wh\mcalx^{(k)},\wh\mcalx)\eem+\frac{\delta+\gamma}{2\delta\gamma}\mbbe\bem\|\mcaly^{(k)}-\wh\mcaly^{(k)}\|_\rmf^2\eem
	\\ & \stackrel{(\ref{ey})(\ref{ed})} \leq \frac{\delta+\gamma}{2\delta\gamma}\frac{\alpha_{\rm r}^2\lambda_{\rm r}\rho_{\rm c}^k}{\|\mcala\|_\rmf^2}\|\mcala*\mcala^\dag*\mcalb\|_\rmf^2+\l(1+\frac{\delta}{\gamma}\r)\rho_{\rm r}\mbbe\bem D_{f,\mcaly^{(k-1)}}(\mcalx^{(k-1)},\wh\mcalx)\eem
	\\ & \leq \frac{\delta+\gamma}{2\delta\gamma}\frac{\alpha_{\rm r}^2\lambda_{\rm r}}{\|\mcala\|_\rmf^2}\|\mcala*\mcala^\dag*\mcalb\|_\rmf^2\l(\rho_{\rm c}^k+\rho_{\rm c}^{k-1}\l(1+\frac{\delta}{\gamma}\r)\rho_{\rm r}\r)\\ &\quad  +\l(1+\frac{\delta}{\gamma}\r)^2\rho_{\rm r}^2\mbbe\bem D_{f,\mcaly^{(k-2)}}(\mcalx^{(k-2)},\wh\mcalx)\eem
	\\ & \leq \cdots 
	\\ & \leq \frac{\delta+\gamma}{2\delta\gamma}\frac{\alpha_{\rm r}^2\lambda_{\rm r}}{\|\mcala\|_\rmf^2}\|\mcala*\mcala^\dag*\mcalb\|_\rmf^2\sum_{i=0}^{k-1}\rho_{\rm c}^{k-i}\l(1+\frac{\delta}{\gamma}\r)^i\rho_{\rm r}^i\\ &\quad  +\l(1+\frac{\delta}{\gamma}\r)^k\rho_{\rm r}^kD_{f,\mcaly^{(0)}}(\mcalx^{(0)},\wh\mcalx).
\end{align*} This completes the proof.
\endproof

\begin{remark}
Let $\rho=\max\{\rho_{\rm c},\rho_{\rm r}\}$. We have $0<\rho<1$. Assume that $\delta>0$ satisfies $0<(1+\delta/\gamma)\rho<1$. We have $$\mbbe\bem D_{f,\mcaly^{(k)}}(\mcalx^{(k)},\wh\mcalx)\eem\leq \l(1+\frac{\delta}{\gamma}\r)^k\rho^k\l(D_{f,\mcaly^{(0)}}(\mcalx^{(0)},\wh\mcalx)+\frac{\delta+\gamma}{2\delta^2}\frac{\alpha_{\rm r}^2\lambda_{\rm r}}{\|\mcala\|_\rmf^2}\|\mcala*\mcala^\dag*\mcalb\|_\rmf^2\r),$$ therefore, 
$$\mbbe\bem\|\mcalx^{(k)}-\wh\mcalx\|_\rmf^2\eem\stackrel{(\ref{gamma})}\leq \l(1+\frac{\delta}{\gamma}\r)^k\rho^k\frac{2}{\gamma}\l(D_{f,\mcaly^{(0)}}(\mcalx^{(0)},\wh\mcalx)+\frac{\delta+\gamma}{2\delta^2}\frac{\alpha_{\rm r}^2\lambda_{\rm r}}{\|\mcala\|_\rmf^2}\|\mcala*\mcala^\dag*\mcalb\|_\rmf^2\r),$$
which shows that the RREK algorithm converges linearly in expectation to $\wh\mcalx$ with the rate $(1+\delta/\gamma)\rho$.	
\end{remark}

\section{Special cases of the proposed algorithm}
\subsection{Tensor randomized extended Kaczmarz (TREK) for tensor least squares} In this subsection, we  consider the following tensor least squares problem \beq\label{tls} \wh\mcalx=\argmin_{\mcalx\in\mbbr^{N_2\times K\times N_3}}\frac{1}{2}\|\mcalx\|_\rmf^2,\quad \mbox{s.t.}\quad \mcala^\top*\mcala*\mcalx=\mcala^\top*\mcalb.\eeq The function $$f(\mcalx)=\frac{1}{2}\|\mcalx\|_\rmf^2$$ is $1$-strongly convex and strongly admissible. We have $$\wh\mcalx=\mcala^\dag*\mcalb,\qquad \nabla f(\mcalx)=\mcalx ,\qquad f^*(\mcalx)=\frac{1}{2}\|\mcalx\|_\rmf^2,\qquad \nabla f^*(\mcalx)=\mcalx,$$ and $$D_{f,\mcalx}(\mcalx,\mcaly)=\frac{1}{2}\|\mcalx-\mcaly\|_\rmf^2.$$ As a consequence, Algorithm 2 becomes Algorithm 3.

\begin{center}
\begin{tabular*}{160mm}{l}
\toprule {\bf Algorithm 3:} The TREK algorithm for solving (\ref{tls})\\ 
\hline \noalign{\smallskip}
\quad {\bf Input}: $\mcala\in\mbbr^{N_1\times N_2\times N_3}$, $\mcalb\in\mbbr^{N_1\times K\times N_3}$, stepsizes $\alpha_{\rm r}$ and $\alpha_{\rm c}$, maximum number of \\ \quad iterations {\tt M}, and tolerance $\tau$.\\\noalign{\smallskip}
\quad {\bf Initialize}: $\mcalz^{(0)}=\mcalb$, $\mcalx^{(0)}\in\ran_K(\mcala^\top)$\\ \noalign{\smallskip}
\quad {\bf for} $k=1,2,\ldots,$ {\tt M} {\bf do}\\ \noalign{\smallskip}
\quad \qquad  Pick $j_k\in[N_2]$ with probability ${\|\mcala_{:,j_k,:}\|^2_\rmf}/{\|\mcala\|_\rmf^2}$\\  \noalign{\smallskip}
\quad \qquad  Set $\mcalz^{(k)}=\mcalz^{(k-1)}-\alpha_{\rm c}\mcala_{:,j_k,:}*\dsp\frac{(\mcala_{:,j_k,:})^\top*\mcalz^{(k-1)}}{\|\mcala_{:,j_k,:}\|_\rmf^2}$\\  \noalign{\smallskip}
\quad \qquad  Pick $i_k\in[N_1]$ with probability ${\|\mcala_{i_k,:,:}\|^2_\rmf}/{\|\mcala\|_\rmf^2}$\\  \noalign{\smallskip}
\quad \qquad Set $\mcalx^{(k)} = \mcalx^{(k-1)}-\alpha_{\rm r}(\mcala_{i_k,:,:})^\top*\dsp \frac{\mcala_{i_k,:,:}*\mcalx^{(k-1)}-\mcalb_{i_k,:,:}+\mcalz^{(k)}_{i_k,:,:}}{\|\mcala_{i_k,:,:}\|^2_\rmf}$ \\  \noalign{\smallskip}
\quad \qquad Stop if $\|\mcalx^{(k)}-\mcalx^{(k-1)}\|_\rmf/\|\mcalx^{(k-1)}\|_\rmf<\tau$ \\  \noalign{\smallskip}
\quad {\bf end}\\
\bottomrule
\end{tabular*}
\end{center}

For any $\mcalx\in\ran_K(\mcala^\top)$, we have $$D_{f,\mcalx}(\mcalx,\wh\mcalx)=\frac{1}{2}\|\mcalx-\wh\mcalx\|_\rmf^2\leq \frac{1}{2\sigma_{\min}^2(\bcirc(\mcala))}\|\mcala*(\mcalx-\wh\mcalx)\|_\rmf^2.$$ Then, the constant $\nu$ in (\ref{nu})  is $\nu=2\sigma_{\min}^2(\bcirc(\mcala))$. In this setting, Theorem \ref{main} reduces to the following result, which implies that the TREK algorithm converges linearly in expectation to the minimum Frobenius norm least squares solution of the tensor system $\mcala*\mcalx=\mcalb$.

\begin{corollary}\label{contrek} Assume that $0<\alpha_{\rm c}<2/\lambda_{\rm c}$ and $0<\alpha_{\rm r}<2/\lambda_{\rm r}$. For any $\delta>0$, the sequence $\{\mcalx^{(k)}\}$  generated by Algorithm 3 satisfies
\begin{align*}\mbbe\bem \|\mcalx^{(k)}-\mcala^\dag*\mcalb\|_\rmf^2\eem &\leq \frac{\delta+1}{\delta}\frac{\alpha_{\rm r}^2\lambda_{\rm r}}{\|\mcala\|_\rmf^2}\|\mcala*\mcala^\dag*\mcalb\|_\rmf^2\sum_{i=0}^{k-1}\rho_{\rm c}^{k-i}\l(1+\delta\r)^i\rho_{\rm r}^i\\ &\quad +(1+\delta)^k\rho_{\rm r}^k\|\mcalx^{(0)}-\mcala^\dag*\mcalb\|_\rmf^2,\end{align*} where $$\rho_{\rm c}=1-\frac{(2\alpha_{\rm c}-\alpha_{\rm c}^2\lambda_{\rm c})\sigma_{\min}^2(\bcirc(\mcala))}{\|\mcala\|_\rmf^2},\quad \rho_{\rm r}=1-\frac{(2\alpha_{\rm r}-\alpha_{\rm r}^2\lambda_{\rm r})\sigma_{\min}^2(\bcirc(\mcala))}{\|\mcala\|_\rmf^2}.$$\end{corollary}

\subsection{RREK for sparse tensor recovery}

In this subsection, we consider the following constrained minimization problem \beq\label{str} \wh\mcalx=\argmin_{\mcalx\in\mbbr^{N_2\times K\times N_3}}\frac{1}{2}\|\mcalx\|_\rmf^2+\lambda\|\mcalx\|_1,\quad \mbox{s.t.}\quad \mcala^\top*\mcala*\mcalx=\mcala^\top*\mcalb.\eeq We mention that minimization problems with the  objective function $$f(\mcalx)=\frac{1}{2}\|\mcalx\|_\rmf^2+\lambda\|\mcalx\|_1$$ have been widely considered; see, e.g., \cite{cai2009conve,lorenz2014linea,schopfer2019linea,chen2021regul} and references therein. Define the soft shrinkage function $S_\lambda(\mcalx)$ componentwise as  $$(S_\lambda(\mcalx))_{i,j,k}= \max\{|\mcalx_{i,j,k}|-\lambda, 0\}\sgn(\mcalx_{i,j,k}),$$ where $\sgn(\cdot)$ is the sign function. We have (see, e.g.,  \cite{beck2017first}) $$\nabla f^*(\mcalx)=S_\lambda(\mcalx).$$ As a consequence, Algorithm 2 becomes Algorithm 4.

\begin{center}
\begin{tabular*}{160mm}{l}
\toprule {\bf Algorithm 4:} The RREK algorithm for solving (\ref{str})\\ 
\hline \noalign{\smallskip}
\quad {\bf Input}: $\mcala\in\mbbr^{N_1\times N_2\times N_3}$, $\mcalb\in\mbbr^{N_1\times K\times N_3}$, stepsizes $\alpha_{\rm r}$ and $\alpha_{\rm c}$, maximum number of \\ \quad iterations {\tt M}, and tolerance $\tau$.\\\noalign{\smallskip}
\quad {\bf Initialize}: $\mcalz^{(0)}=\mcalb$, $\mcaly^{(0)}\in\ran_K(\mcala^\top)$, and $\mcalx^{(0)}=S_\lambda(\mcaly^{(0)})$.\\ \noalign{\smallskip}
\quad {\bf for} $k=1,2,\ldots,$ {\tt M} {\bf do}\\ \noalign{\smallskip}
\quad \qquad  Pick $j_k\in[N_2]$ with probability ${\|\mcala_{:,j_k,:}\|^2_\rmf}/{\|\mcala\|_\rmf^2}$\\  \noalign{\smallskip}
\quad \qquad  Set $\mcalz^{(k)}=\mcalz^{(k-1)}-\alpha_{\rm c}\mcala_{:,j_k,:}*\dsp\frac{(\mcala_{:,j_k,:})^\top*\mcalz^{(k-1)}}{\|\mcala_{:,j_k,:}\|_\rmf^2}$\\  \noalign{\smallskip}
\quad \qquad  Pick $i_k\in[N_1]$ with probability ${\|\mcala_{i_k,:,:}\|^2_\rmf}/{\|\mcala\|_\rmf^2}$\\  \noalign{\smallskip}
\quad \qquad Set $\mcaly^{(k)} = \mcaly^{(k-1)}-\alpha_{\rm r}(\mcala_{i_k,:,:})^\top*\dsp \frac{\mcala_{i_k,:,:}*\mcalx^{(k-1)}-\mcalb_{i_k,:,:}+\mcalz^{(k)}_{i_k,:,:}}{\|\mcala_{i_k,:,:}\|^2_\rmf}$ \\  \noalign{\smallskip}
\quad \qquad Set $\mcalx^{(k)}=S_\lambda(\mcaly^{(k)})$ \\  \noalign{\smallskip}
\quad \qquad Stop if $\|\mcalx^{(k)}-\mcalx^{(k-1)}\|_\rmf/\|\mcalx^{(k-1)}\|_\rmf<\tau$ \\  \noalign{\smallskip}
\quad {\bf end}\\
\bottomrule
\end{tabular*}
\end{center}

Note that the objective function $$f(\mcalx)=\frac{1}{2}\|\mcalx\|_\rmf^2+\lambda\|\mcalx\|_1$$ is $1$-strongly convex and strongly admissible (see \cite[Example 3.7]{chen2021regul} and \cite[Lemma 4.6]{lai2013augme}).  In this setting, Theorem \ref{main} reduces to the following result, which implies the RREK algorithm converges linearly in expectation to the unique solution $\wh\mcalx$. 

\begin{corollary} Let $\nu>0$ be the constant  from (\ref{nu}). Let $\wh\mcalx$ be the solution of (\ref{str}). Assume that $0<\alpha_{\rm c}<2/\lambda_{\rm c}$ and $0<\alpha_{\rm r}<2/\lambda_{\rm r}$. For any $\delta>0$, the sequences $\{\mcalx^{(k)}\}$ and $\{\mcaly^{(k)}\}$ generated by Algorithm 4 satisfy \begin{align*}\mbbe\bem D_{f,\mcaly^{(k)}}(\mcalx^{(k)},\wh\mcalx)\eem &\leq \frac{\delta+1}{2\delta}\frac{\alpha_{\rm r}^2\lambda_{\rm r}}{\|\mcala\|_\rmf^2}\|\mcala*\mcala^\dag*\mcalb\|_\rmf^2\sum_{i=0}^{k-1}\rho_{\rm c}^{k-i}\l(1+\delta\r)^i\rho_{\rm r}^i\\ &\quad +\l(1+\delta\r)^k\rho_{\rm r}^kD_{f,\mcaly^{(0)}}(\mcalx^{(0)},\wh\mcalx),\end{align*} where $$\rho_{\rm c}=1-\frac{(2\alpha_{\rm c}-\alpha_{\rm c}^2\lambda_{\rm c})\sigma_{\min}^2(\bcirc(\mcala))}{\|\mcala\|_\rmf^2},\quad \rho_{\rm r}=1-\frac{(2\alpha_{\rm r}-\alpha_{\rm r}^2\lambda_{\rm r})\nu}{2\|\mcala\|_\rmf^2}.$$\end{corollary}

\section{Numerical experiments}

 In this section, we compare the performance of the proposed RREK algorithm against the RRK algorithm \cite{chen2021regul} in two synthetic problems, including tensor least squares and sparse tensor recovery. The main purpose is to illustrate our theoretical results via simple examples. All experiments are performed using MATLAB R2020b on a laptop with 2.7 GHz Quad-Core Intel Core i7 processor, 16 GB memory, and Mac operating system. The tensor t-product toolbox \cite{lu2018tenso} is used in our computations. In all experiments, the reported results are the average of 10 independent trials. 

\subsection{Tensor least squares} We compare the performance of the TREK algorithm (Algorithm 2) and the TRK algorithm (see section 2.4) for solving the tensor least squares problem (\ref{tls}). 

\begin{figure}[htb]
\centerline{\epsfig{figure=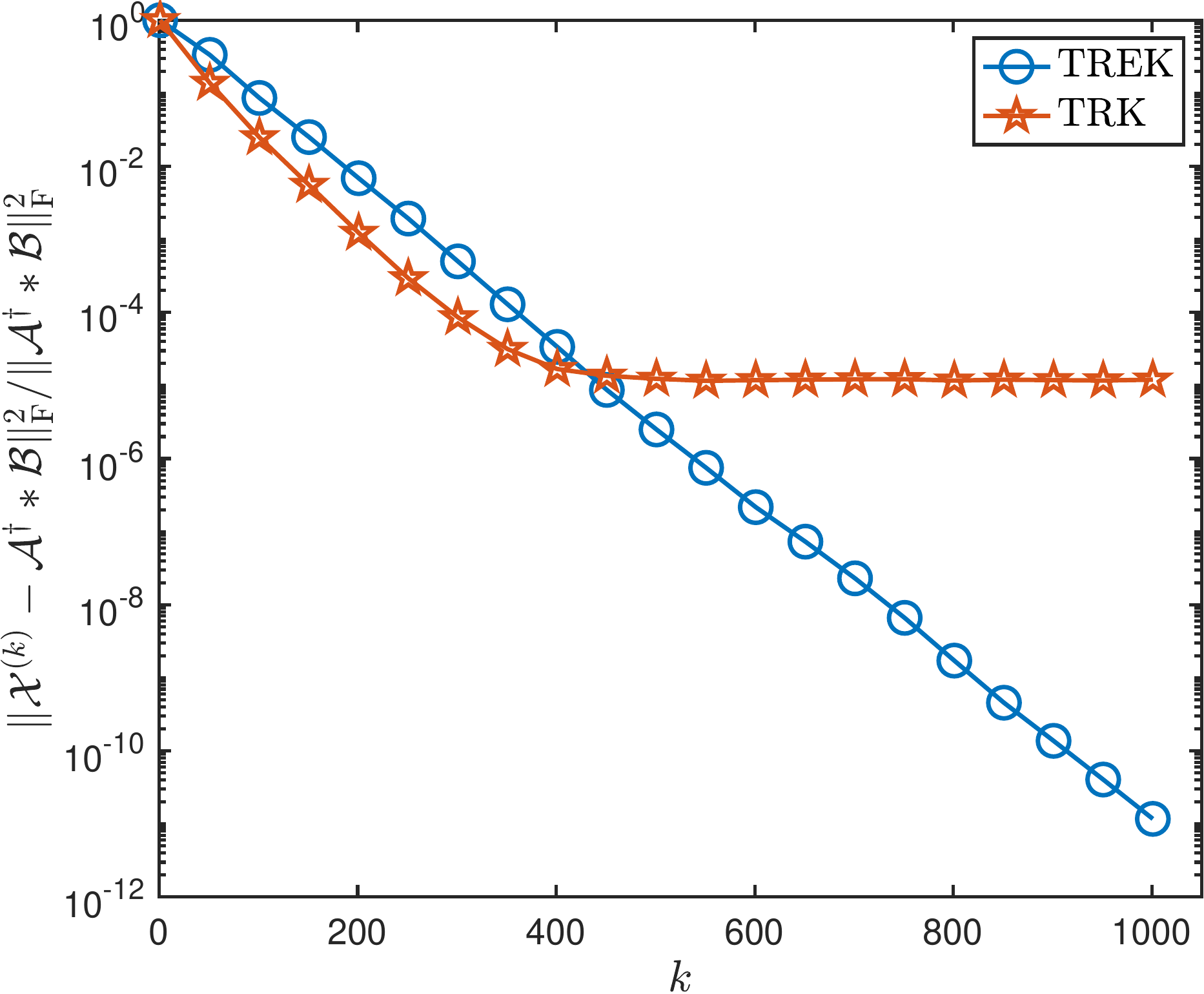,height=3in}}
\caption{The relative error versus the number of iterations of the TREK and TRK algorithms for a tensor least squares problem.}
\label{fig1}
\end{figure}

In our experiment, we generate the sensing tensor $\mcala$ and the acquired measurement tensor $\mcalb$ as follows: $$\mcala={\tt randn}(N_1,N_2,N_3),\quad \mcalb =\mcala*{\tt randn}(N_2,K,N_3)+{\tt randn}(N_1,K,N_3)/10.$$ We set $N_1=100$, $N_2 = 20$, $N_3 = 20$, and $K = 20$. The maximum number of iterations {\tt M} is 1000. We use the stepsizes $\alpha_{\rm r}=1.5/\lambda_{\rm r}$ and $\alpha_{\rm c}=1.5/\lambda_{\rm c}$. In Figure \ref{fig1}, we plot the relative error versus the number of iterations. We observe (i) that the TRK algorithm  converges linearly to within a radius of $\mcala^\dag*\mcalb$ and (ii) that the TREK algorithm converges linearly to $\mcala^\dag*\mcalb$. 

\subsection{Sparse tensor recovery}

We compare the performance of the RREK algorithm (Algorithm 3) and the RRK algorithm \cite[Algorithm 3.2]{chen2021regul} for solving the sparse tensor recovery problem (\ref{str}). 

\begin{figure}[htb]
\centerline{\epsfig{figure=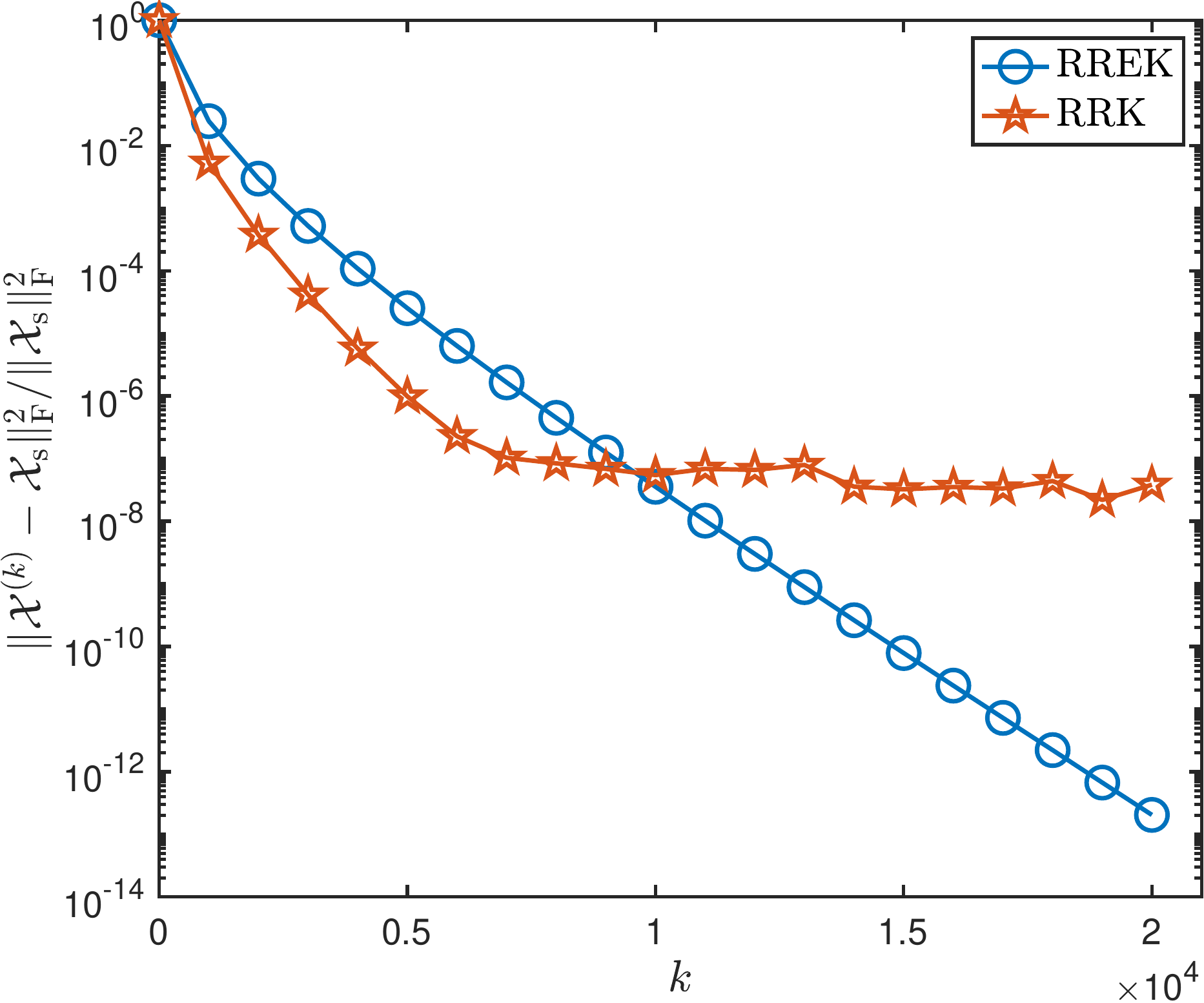,height=3in}}
\caption{The relative error versus the number of iterations of the RREK and RRK algorithms for a sparse tensor recovery problem.}
\label{fig2}
\end{figure}

In our experiment, we generate the sensing tensor $\mcala$, the ground truth tensor $\mcalx_{\rm s}$, and the acquired measurement tensor $\mcalb$ as follows: $$\mcala = {\tt randn}(N_1,N_2,N_3), \quad \mcala(N_1-9:N_1-5,:,:) = \mcala(N_1-4:N_1,:,:),$$ $$\mcalx_{\rm s} = {\tt randn}(N_2,K,N_3),\quad \mcalx_{\rm s} = \mcalx_{\rm s}.*(\mcalx_{\rm s}>=2.33),$$
$$ {\mbf Z} = {\tt null(bcirc}(\mcala)^\top),\quad q = {\tt size}(\mbf Z,2),\quad 
\mcalb = \mcala*\mcalx_{\rm s} + {\tt fold}(\mbf Z*{\tt randn}(q,K))/10.$$ 
We set $N_1=100$, $N_2 = 200$, $N_3 = 10$, and $K = 20$. The maximum number of iterations {\tt M} is 20000. We use the stepsizes $\alpha_{\rm r}=1.5/\lambda_{\rm r}$ and $\alpha_{\rm c}=1.5/\lambda_{\rm c}$. The ground truth $\mcalx_{\rm s}$ is a sparse nonnegative tensor with approximately density $0.01$ and smallest nonzero entry 2.33. In Figure \ref{fig2}, we plot the relative error versus the number of iterations. We observe (i) that the RRK algorithm  converges linearly to within a radius of the ground truth $\mcalx_{\rm s}$ and (ii) that the RREK algorithm converges linearly to the ground truth $\mcalx_{\rm s}$.

\section{Concluding remarks} We have proposed a randomized regularized extended Kaczmarz algorithm for solving the tensor recovery problem (\ref{mp}). At each step, only one horizontal slice and one lateral slice of the sensing tensor $\mcala$ are used. Linear convergence of the proposed algorithm is proved under certain assumptions. Numerical experiments on a tensor least squares problem and a sparse tensor recovery problem confirm the theoretical results. In the future, we will use the existing acceleration strategies such as those in \cite{needell2015rando,du2020rando,bai2021greed,wu2021two} to further improve the efficiency when applied to large scale problems. Moreover, we note that the tensor recovery model (\ref{mp}) can be used as a variable selection procedure and we are studying its performance compared with the elastic net \cite{zou2005regul} and the lasso \cite{tibshirani1996regre}.

\section*{Acknowledgments}
 This work was supported by the National Natural Science Foundation of China (No.12171403 and No.11771364), the Natural Science Foundation of Fujian Province of China (No.2020J01030), and the Fundamental Research Funds for the Central Universities (No.20720210032).	

{\small 

}

\end{document}